\documentclass[review]{elsarticle}

\usepackage{lineno,hyperref}
\modulolinenumbers[1]
\journal{arXiv}

\usepackage{amsmath}
\usepackage{amssymb}
\usepackage[algoruled,boxed,lined]{algorithm2e}
\makeatletter 
\g@addto@macro{\@algocf@init}{\SetKwInOut{Parameter}{Parameters}} 
\makeatother

\usepackage{graphics}
\usepackage{graphicx}
\usepackage{float}

\usepackage{textcomp}

\usepackage{epsfig}
\usepackage{epstopdf}

\usepackage{color}

\usepackage{amsthm}
\usepackage{url}
\usepackage{longtable}
\usepackage[figuresright]{rotating}
\usepackage{listings}
\usepackage{etoolbox}
\usepackage{sectsty}
\usepackage{bm}
\usepackage{algpseudocode}

\usepackage{soul}
\newcommand{\bb}[1]{\mathbb{#1}}

\usepackage{mathtools}
\DeclarePairedDelimiter\abs{\lvert}{\rvert}%

\usepackage[left=3cm, right=3cm, top=3cm]{geometry}

\usepackage{subcaption}

\begin{document}

\begin{frontmatter}

\title{A full multigrid multilevel (Quasi-) Monte Carlo method for the single-phase subsurface flow with random coefficients}


\author[mymainaddress]{Yang Liu}

\author[mymainaddress,mysecondaryaddress]{Jingfa Li}

\author[mymainaddress]{Shuyu Sun\corref{mycorrespondingauthor}}
\cortext[mycorrespondingauthor]{Corresponding author}
\ead{shuyu.sun@kaust.edu.sa}

\author[mysecondaryaddress]{Bo Yu}

\address[mymainaddress]{Computational Transport Phenomena Laboratory, Division of Physical Science and Engineering, King Abdullah University of Science and Technology, Thuwal 23955-6900, Saudi Arabia}
\address[mysecondaryaddress]{School of Mechanical Engineering, Beijing Key Laboratory of Pipeline Critical Technology and Equipment for Deepwater Oil \& Gas Development, Beijing Institute of Petrochemical Technology, Beijing 102617, China}

\begin{abstract}
The subsurface flow is usually subject to uncertain porous media structures. In most cases, however, we only have partial knowledge about the porous media properties. A common approach is to model the uncertain parameters of porous media as random fields, then the statistical moments (e.g. expectation) of the Quantity of Interest(QoI) can be evaluated by the Monte Carlo method. In this study, we develop a full multigrid-multilevel Monte Carlo (FMG-MLMC) method to speed up the evaluation of random parameters effects on single-phase porous flows. In general, MLMC method applies a series of discretization with increasing resolution and computes the QoI on each of them. The effective variance reduction is the success of the method. We exploit the similar hierarchies of MLMC and multigrid methods and obtain the solution on coarse mesh $Q^c_l$ as a byproduct of the full multigrid solution on fine mesh $Q^f_l$ on each level $l$. In the cases considered in this work, the computational saving due to the coarse mesh samples saving is $20\%$ asymptotically. Besides, a comparison of Monte Carlo and Quasi-Monte Carlo (QMC) methods reveals a smaller estimator variance and a faster convergence rate of the latter approach in this study. \\
\end{abstract}

\begin{keyword}
uncertainty quantification \sep subsurface flow \sep Monte Carlo and quasi-Monte Carlo \sep multi-level method \sep multigrid
\end{keyword}

\end{frontmatter}


\section{Introduction}\label{Sec1}

Flow and transport are the most fundamental phenomena in subsurface porous media associated with various physical processes, e.g., oil and gas flow in petroleum reservoir \cite{terry2015applied}, CO$_2$ sequestration \cite{zhang2015sequentially}, water pollution dispersion \cite{bear2010modeling}, etc. The numerical simulation and analyses of flow and transport in subsurface porous media are highly demanded in practical engineering and mechanism studies. However, the simulation results are always subject to the influence of uncertainties, mainly stemming from the inherent spatial heterogeneity of media properties caused by complex geological processes \cite{boschan2012scale}. It has been widely recognized that in natural subsurface porous media, most properties, such as permeability, porosity, etc., exhibit an uneven spatial distribution. For example, the hydraulic conductivity can span several orders of magnitude in an aquifer or reservoir. How to quantificationally identify the influence of uncertainties of porous media properties on the flow and transport behaviors in subsurface physical processes has been a research hot spot in recent years.

Therefore, the uncertainty quantification is an essential task in the simulation of practical subsurface flows where porous media properties that unknown or partially known are taken as the input parameters. A possible way to deal with uncertainties of subsurface porous media is to treat porous media properties as random fields, then perform the stochastic simulation on the subsurface flow governing equations with random coefficients to evaluate the quality of interest (QoI). Among commonly-used stochastic simulation methods, e.g. Monte Carlo(MC) method, stochastic finite element method, stochastic collocation method, the MC method demonstrates apparent advantages such as it is a non-intrusive approach that only the realization of coefficients is needed while the original model code remains unchanged, and it is more easily to be implemented. In the standard MC method, the computer-generated (pseudo) random points are used, and in many cases, the computational efficiency is always unsatisfied for large-scale problems. The Quasi-Monte Carlo (QMC) method improve the demerit of MC by using deterministic quasi-random points. These points exhibit lower discrepancy and distribute more uniformly in the probability space. Moreover, to reduce the sample variance and further improve the computational efficiency, the multilevel Monte Carlo (MLMC) method was proposed and developed by Heinrich \cite{Heinrich2001MultilevelMC} and Giles \cite{giles2008multilevel}. It applies the control variates technique that a series of discretization is adopted with increasing resolution and computes the QoI on each of them, the success of which lies in the effective variance reduction sequentially.

It should be mentioned that in the particular case of subsurface flow with random coefficients, the problem is further aggravated where very detailed geological models are needed (a large number of cells) for an accurate description of the flow. To further alleviate the computational burden connected to the evaluation of random parameter effects on subsurface flow using the MLMC method, in this study, we exploit the similar hierarchies of MLMC and multigrid methods and proposed a full multigrid multilevel (quasi-) Monte Carlo (FMG-MLQMC) approach. In this proposed method, the solution on coarse mesh $Q_l^c$ can be obtained as a byproduct of the full multigrid solution on fine mesh $Q_l^f$ on each mesh level $l$, instead of directly solving the equations on the coarse mesh as the standard MLMC does. The proposed FMG-MLQMC method saves the computation of the $Q_l^c$. There have been works coupling the multigrid solver with the multilevel framework, see \cite{kumar2017multigrid,robbe2018recycling} for example. However, the FMG-MLQMC method we proposed saves the computational cost without modifying the MLMC framework. We exploit the implementation method for upscaling the random coefficient from fine mesh to neighboring coarse mesh. Although in this study we only focus on the simple single-phase subsurface flow with random coefficients, the proposed approach can be applied and extended naturally to multiphase flow and transfer in porous media and any other flow and transport problems associated with uncertainty effect. 

The rest of the paper is organized as follows: In Section \ref{Sec2}, we give a brief description of the single-phase subsurface flow and then introduce the proposed full multigrid multilevel (quasi-) Monte Carlo method in detail. The methodology on the upscaling method of random coefficients from the fine mesh to the coarse mesh, which preserves the random structure, is presented as well. In Section \ref{Sec3}, we verify the effectiveness (a smaller estimator variance and faster convergence rate) of the presented method by comparing with standard MLMC method in two numerical experiments. Finally, in Section \ref{Sec4}, we report the concluding remarks of this work along with a brief discussion of future directions.

\section{Algorithms}\label{Sec2}
\subsection{Model problem and MLMC method}

In this work, we consider the following elliptic problem, 

\begin{equation} \label{eq:elliptic_pde}
    \left\{
    \begin{aligned} 
        -\nabla \cdot (k(\bm{x}, \omega) \nabla u(\bm{x},\omega)) &= f(\bm{x}) \quad \text{ in } D \\
        u(\bm{x},\omega) &= g(\bm{x}) \quad \text{ on } \Gamma_D\\
        \frac{\partial u(\bm{x},\omega)}{\partial \bm{n}} &= v(\bm{x}) \quad \text{ on } \Gamma_N\\
    \end{aligned}
    \right.
\end{equation}
where $k(\bm{x},\omega)$ is the random, spatial-varying coefficient, $D$ is the computational domain, $\omega$ is a sample from the probability triple $(\Omega, \mathcal{F}, P)$. $\Gamma_D$ and $\Gamma_N$ are Dirichlet and Neumann boundaries respectively. In single-phase flow context, Eq.\ref{eq:elliptic_pde} corresponds the steady-state situation, when $g$ and $v$ prescribe the pressure and velocity of the fluid at the boundary, then the solution $u$ depicts the pressure in the domain $\Omega$. 

\begin{figure}[h!]
    \centering
    \includegraphics[width=0.8\textwidth]{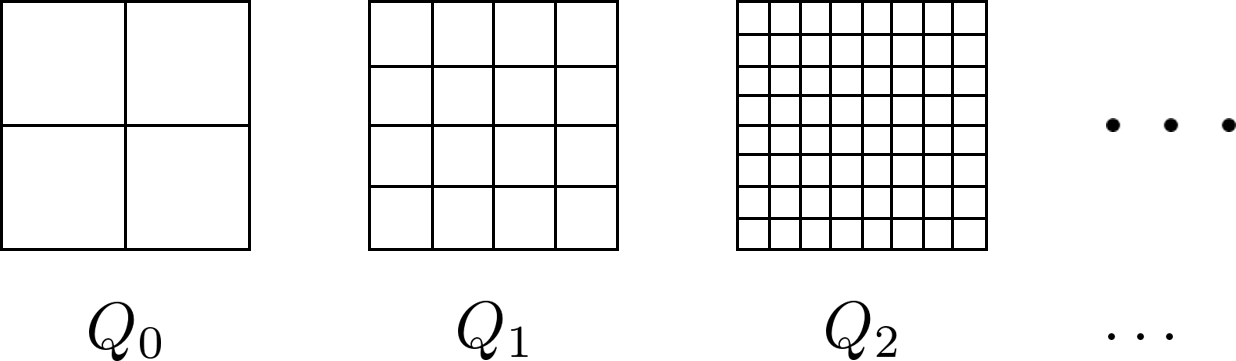}
    \caption{Multilevel Monte Carlo}
    \label{fig:chap2_MLMC}
\end{figure}

In this work, we address the random elliptic problem using the multi-level algorithm. Basically, the MLMC method employs a series of control variates, which are often the discretized models with increasing resolution levels. Here, we associate each level with one mesh with given resolution. The approximations of quantity of interest(QoI) on these levels are denoted as $Q_0,Q_1,\cdots,Q_L$, see Figure \ref{fig:chap2_MLMC}. 

We are interested in the approximation $Q_L$ on the finest level $L$. The MLMC method not only computes the solution on level $L$ itself, but also calculates the solutions on all the preceding meshes. The expectation of such quantity can be expressed by the following telescoping formula,
\begin{equation} \label{eq:Chap2_Telescoping_Formula}
    \bb{E}[Q_L] = \bb{E}[Q_0] + \sum_{l=1}^{L} (\bb{E}[Q_l] - \bb{E}[Q_{l-1}])
\end{equation}

We can approximate each expectation using the Monte Carlo approach as follows,
\begin{equation} \label{eq:Chap2_Telescoping_Formula_approx}
\bb{E}[Q_L]\approx \frac{1}{N_0}\sum_{i=1}^{N_0}Q_0(\omega_{i,0})+\sum_{l=1}^{L}\frac{1}{N_l}[\sum_{i=1}^{N_l}(Q_l(\omega_{i,l})-Q_{l-1}(\omega_{i,l})]
\end{equation}

Here, we associate level $l$ with the $l$-th term in the telescoping formula \eqref{eq:Chap2_Telescoping_Formula}. Notice that on each level $l$, we use the same samples $\omega_{i,l}$ to calculate $Q_l$ and $Q_{l-1}$. Then $Q_l$ and $Q_{l-1}$ are likely to correlate well, and the variance of $Q_l - Q_{l-1}$ will be small, see Eq.\ref{eq:Chap2_variance_small}. 

\begin{equation}\label{eq:Chap2_variance_small}
\begin{split}
V_l = \bb{V}[Q_l - Q_{l-1}] &= \bb{V} [Q_l] + \bb{V} [Q_{l-1}] - 2Cov(Q_l,Q_{l-1})\\
& \ll \bb{V} [Q_l] + \bb{V} [Q_{l-1}]
\end{split}
\end{equation}
where we denote $V_l$ as the variance of the $Q_l - Q_{l-1}$ on level $l$. 

Also notice that on different levels, independent samples are used so that the variance of the multilevel estimator $Q_L$ is the summation of the variance on each level.\\
\begin{equation} \label{eq:chap2_variance_sum}
    \bb{V}[Q_L] = \sum_{l=0}^{L} \bb{V}[Q_l - Q_{l-1}] = \sum_{l=0}^L V_l,
\end{equation}
where we let $Q_{-1} = 0$. If we write $Y_0=Q_0$ and $Y_l=Q_l-Q_{l-1}$, then
\begin{equation*}
\bb{E}[Q_L]=\sum_{l=0}^{L}\bb{E}[Y_l]
\end{equation*}
Let $\hat{Y}_l$ be an unbiased estimator for $E[Y_l]$, 
\begin{equation*}
\begin{aligned}
\hat{Y}_0 &=\frac{1}{N_0}\sum_{i=1}^{N_0}Q_0(\omega_{i,0})\\
\hat{Y}_l &=\frac{1}{N_l}[\sum_{i=1}^{N_l}(Q_l(\omega_{i,l})-Q_{l-1}(\omega_{i,l})] \qquad {l=1,2,3,\cdots,L}
\end{aligned}
\end{equation*}
then the multilevel estimator becomes,
\begin{equation}
\hat{Q}^{ML}_L=\sum_{l=0}^{L}\hat{Y}_l    
\end{equation}

\subsection{MLMC Complexity Theory}
Let $Q$ denote a quantity of interest, and $Q_l$ denote the corresponding numerical approximation on $l$-th mesh. If we assume that the weak error and the level variance decreases exponentially while the cost per sample on each level increases exponentially, there exist positive constants $\alpha$, $\beta$ and $\gamma$ satisfying the following,
\begin{equation}
    \begin{split}
        \left \lvert \bb{E}[Q_l-Q] \right\rvert &= \mathcal{O}(2^{-\alpha l})\\
        \bb{V}[Y_l] &= \mathcal{O}(2^{-\beta l})\\
        C_l &= \mathcal{O}(2^{\gamma l})
    \end{split}
    \label{eq:chap2_assumption}
\end{equation}
where $C_l$ is the cost per sample on level $l$. With the mean square error less than a threshold,

\begin{equation}
\bb{E}[(\sum_{l=0}^{L} \hat{Y}_l - \bb{E}[Q])^2] = \sum_{l=0}^L N_l^{-1} V_l+(\bb{E}[\hat{Q}^{ML}_L-Q])^2 < \epsilon ^{2}\\
\label{eq:chap2_mlmc_computation_goal}
\end{equation}
the total computational cost satisfies,
\[
C=\left\{
    \begin{array}{ll}
         \mathcal{O}(\epsilon^{-2}) & \beta > \gamma\\
         \mathcal{O}(\epsilon^{-2} (\log \epsilon )^2) & \beta = \gamma\\
         \mathcal{O}(\epsilon^{-2-(\gamma-\beta)/\alpha})  & \beta < \gamma\\
    \end{array}
\right.
\]
as $\epsilon \to 0$. 


\subsection{MLMC Algorithm} \label{Chapter:MLMC_Algorithm}
This subsection gives a MLMC algorithm initially proposed by M.Giles\cite{giles2008multilevel}.\\

\begin{algorithm}[H]
\SetAlgoLined
 Start with $L=2$, set the initial number of samples on level 0,1,2\\
 \While{extra samples need to be evaluated}{
  evaluate $Q_l(\omega_{i,l})$ and $Q_{l-1}(\omega_{i,l})$, for $\{i,l:dN_l \neq 0, i = 1,\cdots,dN_l\}$\\
  update estimates for $V_l$, $l=0,\cdots,L$\\
  update optimal $N_l$, compute the number of extra samples $dN_l$\\
  \eIf{$|\bb{E}[Q_L-Q]| \approx \frac{|\bb{E}[Q_L-Q_{L-1}]|}{(2^{\alpha}-1)}<\frac{\epsilon}{\sqrt{2}}$}{
   \textbf{break}
   }{
   $L=L+1$ and initialize $N_L$\;
  }
 }
 \caption{MLMC Algorithm}
 \label{alg:MLMC}
\end{algorithm}

In Algorithm \ref{alg:MLMC}, the variances $V_l$ are approximated by the sample variances on the run. The weak error $\abs{\bb{E}[Q_L-Q]}$ is approximated by Richardson extrapolation $\frac{\abs{\bb{E}[Q_L-Q_{L-1}]}}{(2^{\alpha}-1)}$. 

And here we consider a equal split of the estimator variance and the approximation error, i.e.,
\begin{align}
    \sum_{l=0}^L N_l^{-1} V_l &< \epsilon^2/2\\
    \bb({E}[\hat{Q}^{ML}_L - Q])^2 &< \epsilon^2/2. \label{eq:Chap_2_variance_constraint}
\end{align}
However, it is possible to determine the split factor in an optimal way\cite{Collier2015}. 

The optimal number of samples $N_l$ can be obtained by solving a constrained optimization : minimizing the total computational cost subject to the constraint Eq.\ref{eq:Chap_2_variance_constraint}. 

MLMC algorithm will work under the following three conditions. 

\paragraph{Convergence} The sequence $Q_0, Q_1, \cdots, Q_L, \cdots$ converges. Otherwise, the telescoping equation \eqref{eq:Chap2_Telescoping_Formula} does not yield a converging result. 
\paragraph{Correlation} $Q_l$ and $Q_{l-1}$ are estimated using the same underlying random sample $\omega_{i,l}$ in equation \eqref{eq:Chap2_Telescoping_Formula_approx}, and are thus well correlated. In this case, the estimator variance is significantly reduced. 
\paragraph{Consistency} The telescoping sum \eqref{eq:Chap2_Telescoping_Formula} introduces no bias error. Notice that in the telescoping equation, the term $Q_{l-1}$ for $l=1,\cdots,L$ appears twice. However, the two $Q_{l-1}$ may be evaluated differently. If we denote the $l$-th term in the telescoping equation $Q_{l}$ and $Q_{l-1}$ by $Q_{l}^{f}$ and $Q_{l}^{c}$, respectively, which denote the fine mesh solution and coarse mesh solution on level $l$, then the condition $\bb{E}[Q_{l-1}^f] = \bb{E}[Q_l^c]$ needs to be satisfied in order to introduce no bias error in equation \eqref{eq:Chap2_Telescoping_Formula}. The expectation of fine mesh solution on level $l-1$ should be the same as that of the coarse mesh solution on level $l$.

\subsection{MLQMC Algorithm}

\begin{figure}[h!]
    \centering
    \begin{subfigure}{.33\textwidth}
        \centering
        \includegraphics[width=0.97\textwidth]{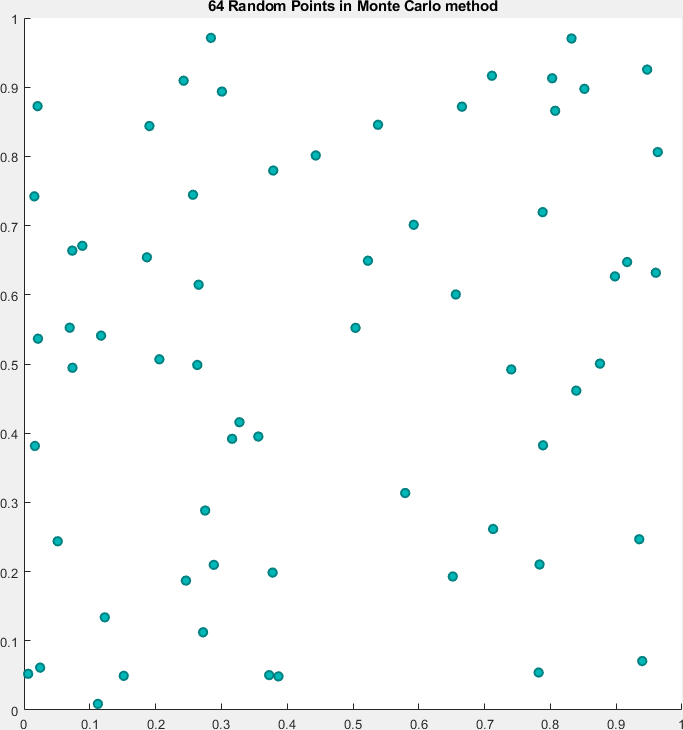}
        \caption{Monte Carlo}
    \end{subfigure}
    \begin{subfigure}{0.33\textwidth}
        \centering
        \includegraphics[width=0.97\textwidth]{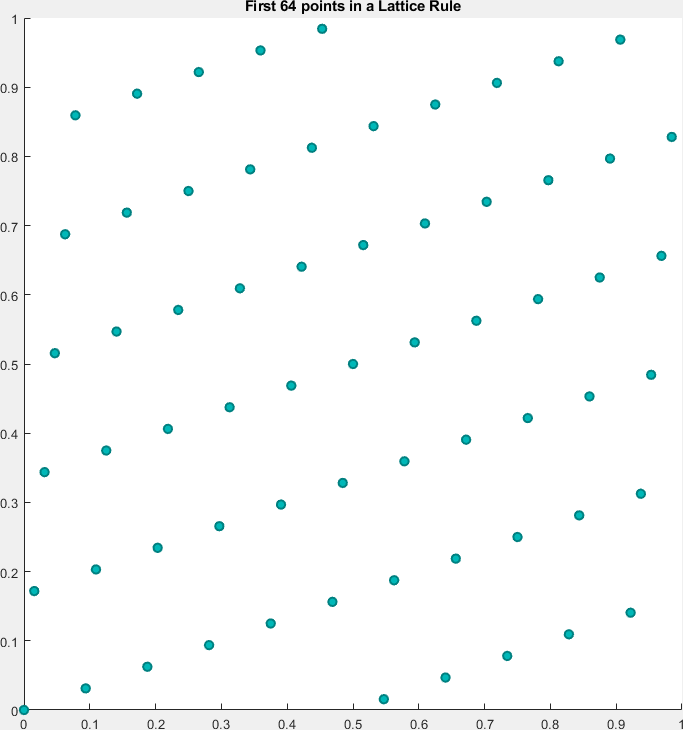}
        \caption{Lattice rule}
    \end{subfigure}\hfill
    \begin{subfigure}{0.33\textwidth}
        \centering
        \includegraphics[width=0.97\textwidth]{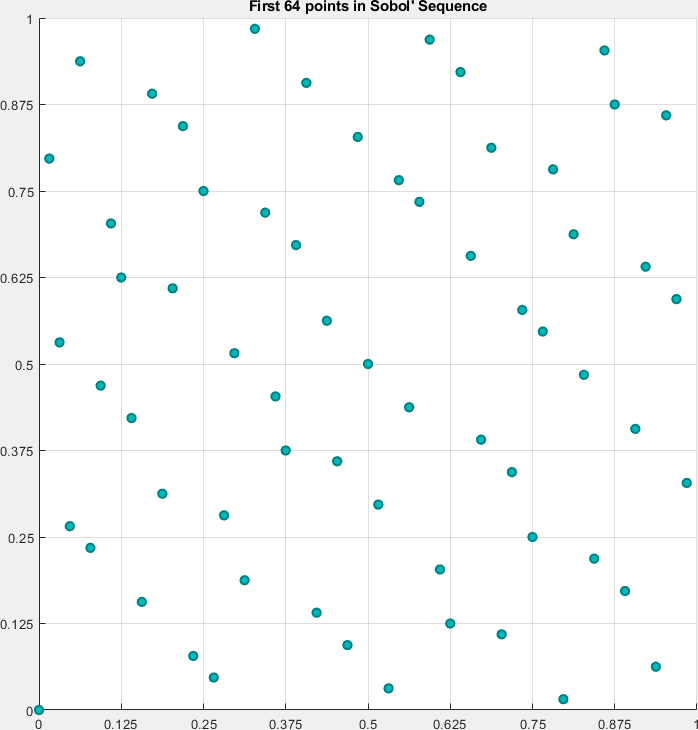}
        \caption{Sobol' sequence}
    \end{subfigure}
    \caption{An example of Monte Carlo points, Lattice rule and Sobol' sequence in $[0,1]\times[0,1]$ domain.}
    \label{fig:MC_Lattice_Sobol}
\end{figure}

The QMC method can solve integration problems as well. In contrast to the MC method, the QMC method replaces random points with deterministic points. Figure \ref{fig:MC_Lattice_Sobol} gives an example of Monte Carlo points, lattice rule, and Sobol' sequence. 

The QMC approximation is given by,
\[
\mathcal{I}_{QMC}^{N}=\frac{1}{N} \sum_{i=1}^{N} Q(t_i) \approx \bb{E}[Q].
\]

Notice that the points $\{t_i\}_{i=1}^{N}$ are deterministic. However, deterministic points yield biased estimates. In this work, we use the randomized QMC: random shift and digital scramble for lattice rule and Sobol' sequence respectively. The interested readers may refer to \cite{dick2013high} for the above mentioned two randomization techniques. 

There have been numerous studies combining the MLMC and QMC methods\cite{giles2009multilevel,giles2016combining,kuo2017multilevel}. 
We follow the multilevel quasi-Monte Carlo (MLQMC) settings from these works and list the algorithm here.

\begin{algorithm}[H]
\SetAlgoLined
Start with $L=2$, set initial number of samples $N_0$ on level 0,1,2\\
\While{extra samples need to be evaluated}{
	evaluate $Q_l(\omega_{i,l})$ and $Q_{l-1}(\omega_{i,l})$, for $\{i,l:dN_l \neq 0, i = 1,\cdots,dN_l\}$\\
	update estimates for $V_l$, $l=0,\cdots,L$ and compute $\bb{V}[Q]$\\
	\eIf{$\bb{V}[Q] > \epsilon^2/2$}{
		select level $l$ such that $l=\text{argmax} \frac{\bb{V}[Y_l]}{N_lC_l}$ and double $N_l$}{
		\eIf{$\bb{E}[|Q_L-Q|]\approx\frac{\abs{\bb{E}[Q_L-Q_{L-1}]}}{(2^{\alpha}-1)}<\frac{\epsilon}{\sqrt{2}}$}{
			\textbf{break}
		}{
			$L=L+1$ and initialize $N_L$ 
		}
	}
}
\caption{MLQMC Algorithm}
\end{algorithm}

\subsection{Multigrid}

\begin{figure}[H]
    \centering
    \includegraphics[width=0.80\textwidth]{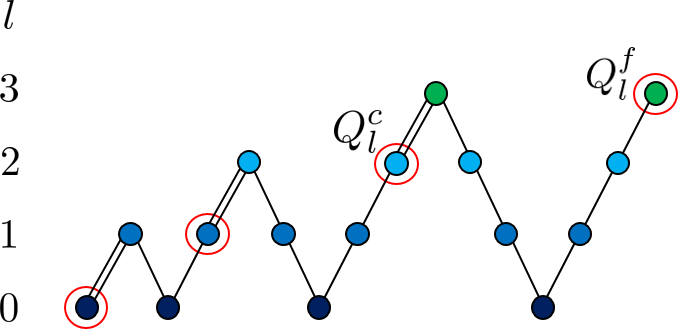} 
    \caption{An illustration of Full-Multigrid-Multilevel Monte Carlo method.}
    \label{fig:FMG-MLMC}
\end{figure}

The multigrid method was originally introduced to solve elliptic boundary-value problems efficiently. It has since been developed to solve either linear or non-linear systems. Multigrid methods compute the solution on a sequence of grids. Figure. \ref{fig:FMG-MLMC} gives an illustration of the full multigrid scheme.

We observe that, when the full multigrid solver is applied to the MLMC problem, based on the same level hierarchies, the solution on the coarse mesh $Q^c_l$ can be obtained as a byproduct of the multigrid solution on fine mesh $Q^f_l$ on each level $l$. Thus, in our proposed FMG-MLMC method, we have saved the computation for $Q^c_l$. Also notice that at the red circles in Figure. \ref{fig:FMG-MLMC} the solutions are exact, since they are the end point of each V-cycle. 

Recall that $Q_l(\omega_{i,l})$ and $Q_{l-1}(\omega_{i,l})$ are evaluated by the same underlying random sample. In level $l$, we denote $K_l^f$ and $K_l^c$ as the coefficients of the fine and coarse models, respectively. Also recall the consistency condition, since we use the same numerical solver for $Q_l^c$ and $Q_{l-1}^f$, we only require that $K_l^c$ and $K_{l-1}^f$ follow the same distribution. 

\begin{figure}[H]
    \centering
    \includegraphics[width=0.6\textwidth]{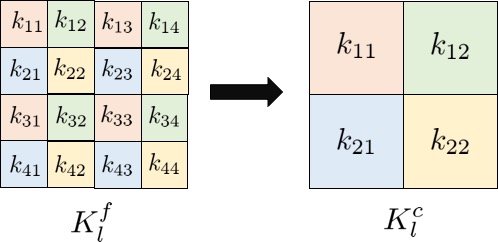}
    \caption{Coarsening}
    \label{fig:chap3_coarsening}
\end{figure}

$K_l^f$ can be generated by the matrix decomposition method, KL-expansion method or other random field generation methods, see \cite{Liu2019} for example. A way to generate $K_l^c$ is to coarsen $K_l^f$. In order to prevent bias error, $K_l^c$ should satisfy the same distribution law as $K_l^f$. Figure \ref{fig:chap3_coarsening} shows a way of coarsening. 

In this scheme, the value of the coefficient in the coarse grid is selected to be the corresponding block in the fine grid. We denote the blocks in fine grid and coarse grid by $k^f_{i,j}$ and $k^c_{I,J}$ respectively, then we have the following,

\begin{equation*}
    k^c_{I,J} = k^f_{2I-1, 2J-1}.
\end{equation*}

Here we give an algorithm to describe our proposed FMG-MLMC method. 

\begin{algorithm}[H]
\SetAlgoLined
 Start with $L=2$, set the initial number of samples on level 0,1,2\\
 \While{extra samples need to be evaluated}{
  coarsen $K_l^f$ to obtain $K_l^c$\\
  use the multigrid solver to compute realizations $Q_l(\omega_{i,l})$, and obtain $Q_{l-1}(\omega_{i,l})$ as a byproduct, for $\{i,l:dN_l \neq 0, i = 1,\cdots,dN_l\}$\\
  update estimates for $V_l$, $l=0,\cdots,L$\\
  update optimal $N_l$, compute number of extra samples $dN_l$\\
  \eIf{$|\bb{E}[Q_L-Q]| \approx \frac{|\bb{E}[Q_L-Q_{L-1}]|}{(2^{\alpha}-1)}<\frac{\epsilon}{\sqrt{2}}$}{
   \textbf{break}
   }{
   $L=L+1$ and initialize $N_L$\;
  }
 }
 \caption{FMG-MLMC Algorithm}
\end{algorithm}

Notice that in our scheme, no correlation is introduced into each levels. Thus the equation \eqref{eq:chap2_variance_sum} still holds true. 

\section{Numerical Validation}\label{Sec3}
\subsection{Problem Statement}\label{Sec3.1}
Recall the elliptic problem \ref{eq:elliptic_pde} in Section 2, now we focus the physical domain $\Omega=[0,1]^2$. In this work, we consider cases in two different boundary conditions and quantities of interest, as listed in Table \ref{tab:Chap4_cases}. Case I is of the Dirichlet boundary type, with pointwise output quantity. Case II is of the mixed Dirichlet-Neumann boundary condition, whose output is the outflow at the east boundary.

\begin{table}[h]
\centering
\caption{Case Settings}
\begin{tabular}{ccc}\hline
Case & Boundary Condition & QoI \\ \hline
1 & $u\mid_{\partial W} = 100$, $u\mid_{\partial E} = 0$, $u\mid_{\partial N} = 50$, $u\mid_{\partial S} = 10$ & $u(0.5,0.5)$ \\
2 & $u\mid_{\partial W} = 100$, $u\mid_{\partial E} = 0$, $\frac{\partial u}{\partial n}\mid_{\partial N}=0$, $\frac{\partial u}{\partial n}\mid_{\partial S}=0$ & $\int_{\partial E} -k\nabla udx$ \\\hline
\end{tabular}
\label{tab:Chap4_cases}
\end{table}

\null
\noindent (1) \textbf{Discretization}\\
The governing equation (\ref{eq:elliptic_pde}) is discretized by the finite-volume method on rectangular grids. On level $l$ the degree of freedom is $2^{l+2} \times 2^{l+2}$.\\
(2) \noindent\textbf{Random Fields}\\
We choose the Mat\'ern covariance function 
\begin{equation} \label{eq:Matern Covariance}
C_\nu (d)=\sigma^2 \frac{2^{1-\nu}}{\Gamma(\nu)} (\sqrt{2\nu} \frac{d}{\lambda})^{\nu} K_\nu (\sqrt{2\nu}\frac{d}{\lambda}),
\end{equation}
where $d$ is the Euclidean distance of two points, $\lambda$ is the correlation length, and $\nu$ controls the smoothness of the field.

The parameters of Mat\'ern covariance for four different random fields are given in Table \ref{table:Random Field Parameters Setting}.
\begin{table}[h] 
\centering
\caption{Random Field Parameters Settings}
\label{table:Random Field Parameters Setting}
\begin{tabular}{cc}\hline
Random Field & Parameters \\ \hline
1 & $\nu = 0.5, \lambda = 0.5, \sigma^2 = 1$ \\
2 & $\nu = 0.5, \lambda = 1, \sigma^2 = 1$\\
3 & $\nu = 1, \lambda = 0.5, \sigma^2 = 1$\\
4 & $\nu = 1, \lambda = 1, \sigma^2 = 1$\\\hline
\end{tabular}
\end{table} 

The random fields are generated using the KL-expansion method. The truncation term is determined when $99\%$ of the variability is captured, meaning that
\[
\frac{\sum_{i=1}^{N_{KL}} \theta_i}{\sum_{i=1}^{\infty} \theta_i} = 99\%,
\]
where the summation of all eigenvalues satisfies the following
\[
\sum_{i=1}^{\infty} \theta_i = \sigma^2 meas(\Omega) = \sigma^2 \int_{\Omega} dx,
\]
Readers of interest can see \cite{ernst2009efficient} for examples. $\Omega$ is the random field region. \\
(3) \noindent\textbf{QMC Method}\\
In QMC method, the Lattice rule points are generated using the software from \cite{kuo2016application}, Sobol' matrices from \cite{joe2008constructing} are used to generate Sobol' sequences. In both cases, 24 randomizations are applied. The confidence intervals are obtained by 10 sets of randomization.

\subsection{Numerical results}\label{Sec3.2}

In this subsection we present the numerical results of the two cases (Tab. \ref{tab:Chap4_cases}) with four random field settings (Tab. \ref{table:Random Field Parameters Setting}). In each figure, the first and second row corresponds the cases $\nu = 0.5$ and $\nu = 1.0$ respectively, while the first column and second column corresponds the cases $\lambda = 0.5$ and $\lambda = 1.0$. We will first give the results of the first case in the following. 

\begin{figure}[H]
    \begin{center} 
    \hfill
    \begin{minipage}{0.33\textwidth}
        \centering
        \includegraphics[width=\textwidth]{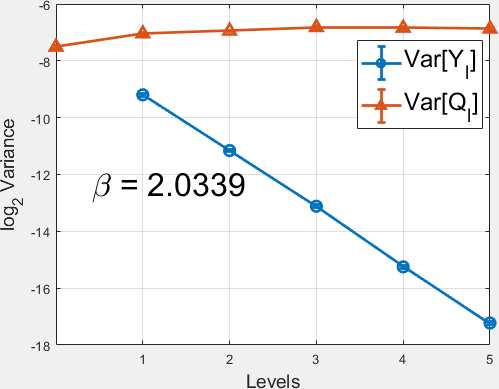} 
    \end{minipage}
    \begin{minipage}{0.33\textwidth}
        \centering
        \includegraphics[width=\textwidth]{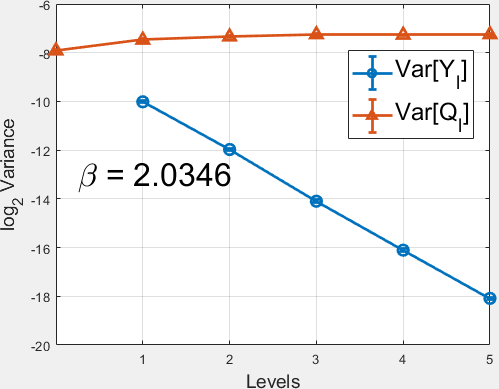} 
    \end{minipage}\hfill
    \null \\
    \begin{minipage}{0.33\textwidth}
        \centering
        \includegraphics[width=\textwidth]{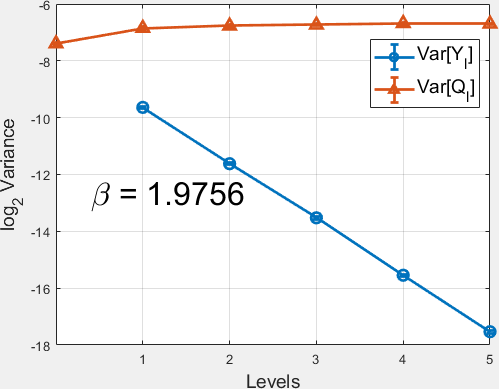} 
    \end{minipage}
    \begin{minipage}{0.33\textwidth}
        \centering
        \includegraphics[width=\textwidth]{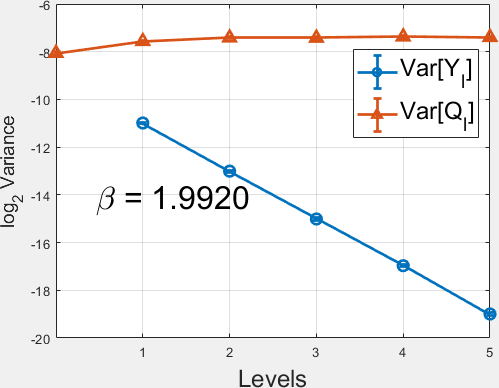}
    \end{minipage}\hfill
    \end{center}
    \caption{Variance of $Q_l$ and $Y_l$ for 4 random fields}
    \label{fig:mlmc_1_var}
\end{figure}

\begin{figure}[H]
    \begin{center}
    \hfill
    \begin{minipage}{0.33\textwidth}
        \centering
        \includegraphics[width=\textwidth]{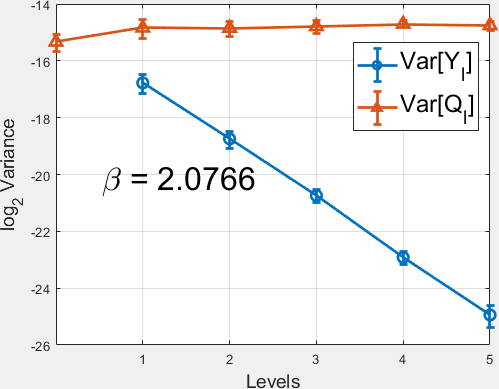} 
    \end{minipage}
    \begin{minipage}{0.33\textwidth}
        \centering
        \includegraphics[width=\textwidth]{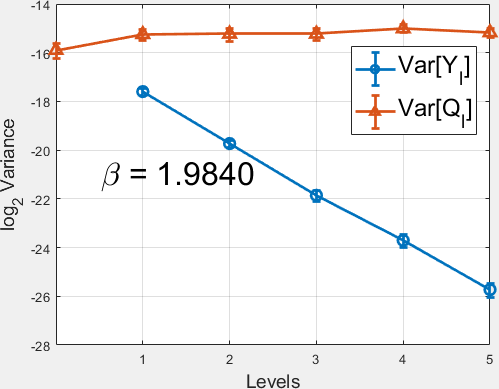} 
    \end{minipage}\hfill
    \null \\
    \begin{minipage}{0.33\textwidth}
        \centering
        \includegraphics[width=\textwidth]{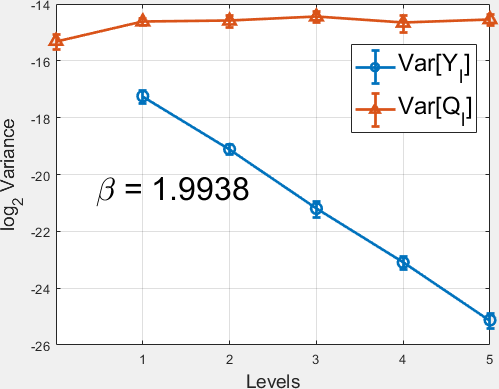} 
    \end{minipage}
    \begin{minipage}{0.33\textwidth}
        \centering
        \includegraphics[width=\textwidth]{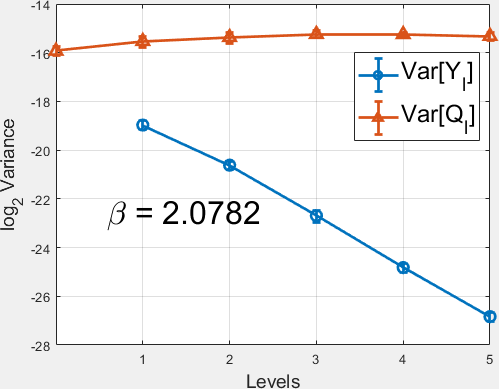}
    \end{minipage}\hfill
    \end{center}
    \caption{Variance of $Q_l$ and $Y_l$ for 4 random fields}
    \label{fig:mlqmc_lattice_1_var}
\end{figure}

\begin{figure}[H]
    \begin{center}
    \hfill
    \begin{minipage}{0.33\textwidth}
        \centering
        \includegraphics[width=\textwidth]{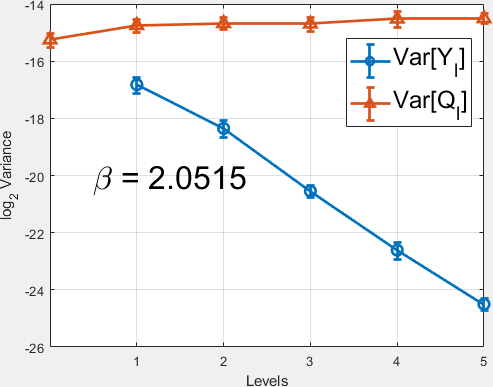} 
    \end{minipage}
    \begin{minipage}{0.33\textwidth}
        \centering
        \includegraphics[width=\textwidth]{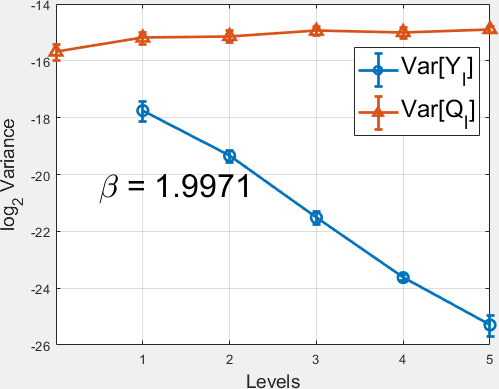} 
    \end{minipage}\hfill
    \null \\
    \begin{minipage}{0.33\textwidth}
        \centering
        \includegraphics[width=\textwidth]{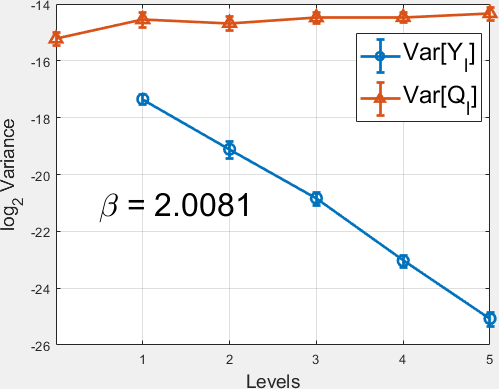} 
    \end{minipage}
    \begin{minipage}{0.33\textwidth}
        \centering
        \includegraphics[width=\textwidth]{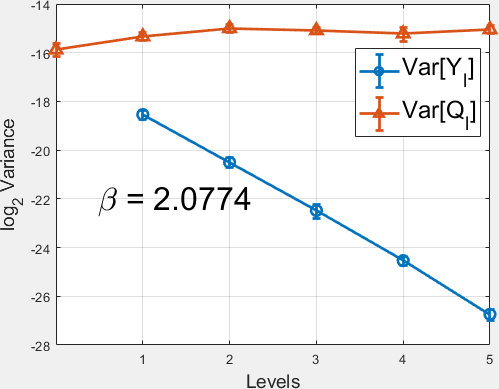}
    \end{minipage}\hfill
    \end{center}
    \caption{Variance of $Q_l$ and $Y_l$ for 4 random fields}
    \label{fig:mlqmc_sobol_1_var}
\end{figure}

The simulation starts by estimating the asymptotic rates $\alpha, \beta, \gamma$ in the assumptions \eqref{eq:chap2_assumption}. Figs. \ref{fig:mlmc_1_var}, \ref{fig:mlqmc_lattice_1_var} and \ref{fig:mlqmc_sobol_1_var} plot the variance of the QoI $Q_l$ and $Y_l = Q_l - Q_{l-1}$ against level $l$. By comparison, the QMC method reduces in the variance not in the asymptotic variance convergence rate, but in the y-axis offsets. 

For the results of MLQMC-Lattice, the variance test (Fig. \ref{fig:mlqmc_lattice_1_var}) is presented. 

For the results of MLQMC-Sobol', the variance (Fig. \ref{fig:mlqmc_sobol_1_var}) is presented. 

\begin{figure}[H]
    \begin{center}
    \hfill
    \begin{minipage}{0.33\textwidth}
        \centering
        \includegraphics[width=\textwidth]{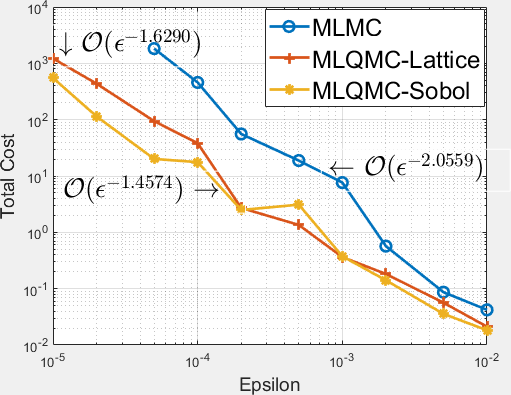} 
    \end{minipage}
    \begin{minipage}{0.33\textwidth}
        \centering
        \includegraphics[width=\textwidth]{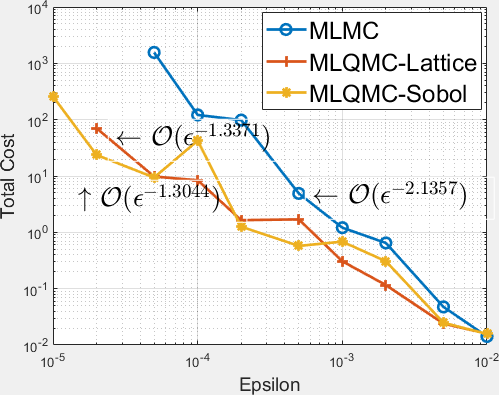} 
    \end{minipage}\hfill
    \null \\
    \begin{minipage}{0.33\textwidth}
        \centering
        \includegraphics[width=\textwidth]{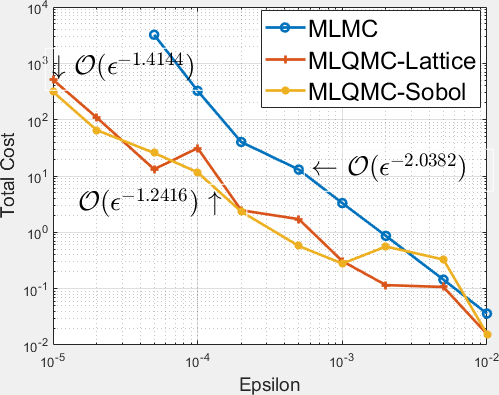} 
    \end{minipage}
    \begin{minipage}{0.33\textwidth}
        \centering
        \includegraphics[width=\textwidth]{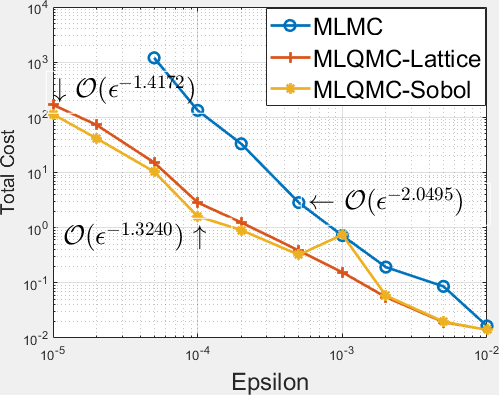}
    \end{minipage}\hfill
    \end{center}
    \caption{Computational complexity of the three methods; the asymptotic rates are marked on the plot}
    \label{fig:case1_comparison}
\end{figure}

The QMC method does not affect the expectation value nor the computational cost. Here we skip the comparison of $\alpha$ and $\gamma$, but present the final computational cost of the three methods, given $\epsilon$. The results in 4 random fields are presented in Fig. \ref{fig:case1_comparison}. 

For case II, we test the QMC convergence rate and then present the computational complexity. 

\begin{figure}[H]
    \begin{center}
    \hfill
    \begin{minipage}{0.33\textwidth}
        \centering
        \includegraphics[width=\textwidth]{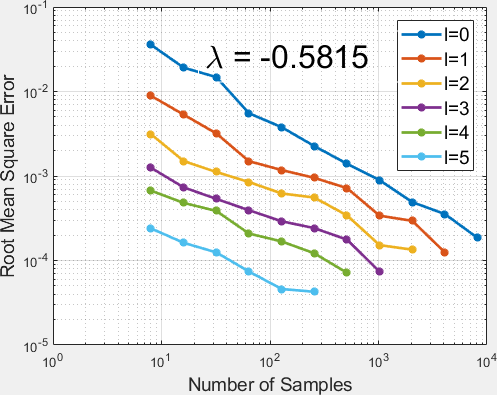} 
    \end{minipage}
    \begin{minipage}{0.33\textwidth}
        \centering
        \includegraphics[width=\textwidth]{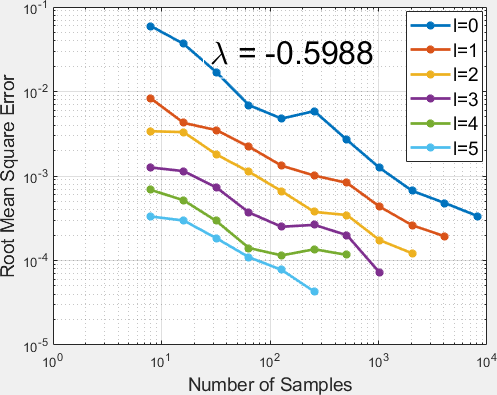} 
    \end{minipage}\hfill
    \null \\
    \begin{minipage}{0.33\textwidth}
        \centering
        \includegraphics[width=\textwidth]{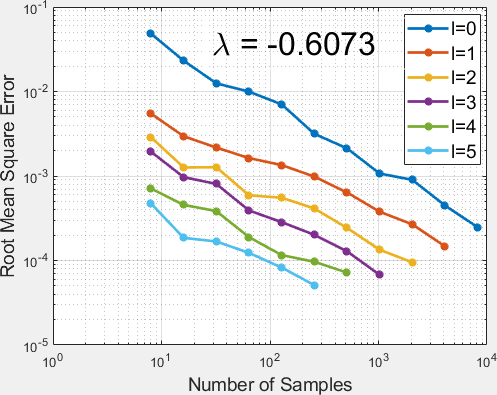} 
    \end{minipage}
    \begin{minipage}{0.33\textwidth}
        \centering
        \includegraphics[width=\textwidth]{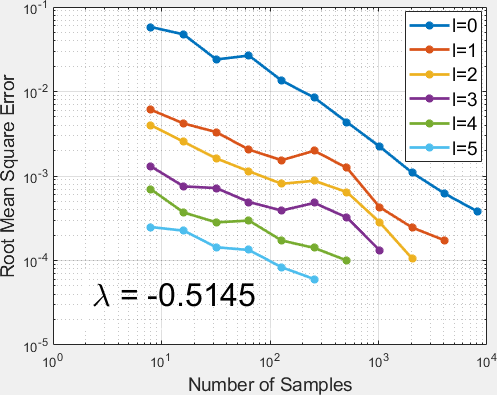}
    \end{minipage}\hfill
    \end{center}
    \caption{Variance test of MLQMC-Lattice for all levels. The variance of the estimator is plotted as a function of the number of samples $N_l$. The convergence rates $\lambda$ are marked on the plot}
    \label{fig:mlqmc_lattice_2_test}
\end{figure}

Next the results of MLQMC-Lattice, the convergence test (Fig. \ref{fig:mlqmc_lattice_2_test}) variance of level estimator $Y_l$ against $N_l$ in each level. The offset between the lines reveals the variance decreases with the levels. The comparison shows Sobol' sequence's advantage in decreasing variance. The random parameter settings have small impact on the variances in this case. 

\begin{figure}[H]
    \begin{center}
    \hfill
    \begin{minipage}{0.33\textwidth}
        \centering
        \includegraphics[width=\textwidth]{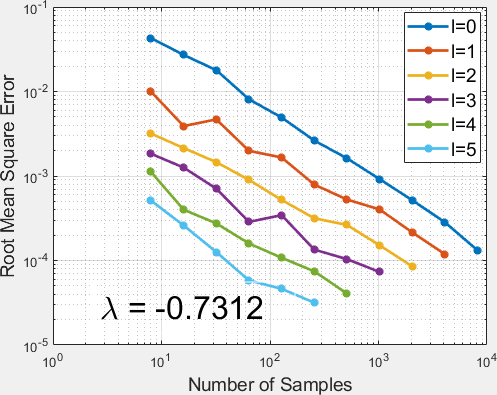} 
    \end{minipage}
    \begin{minipage}{0.33\textwidth}
        \centering
        \includegraphics[width=\textwidth]{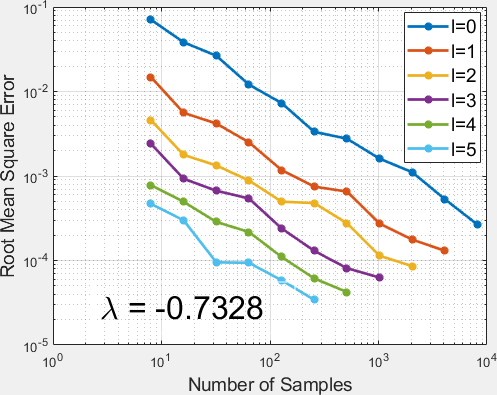} 
    \end{minipage}\hfill
    \null \\
    \begin{minipage}{0.33\textwidth}
        \centering
        \includegraphics[width=\textwidth]{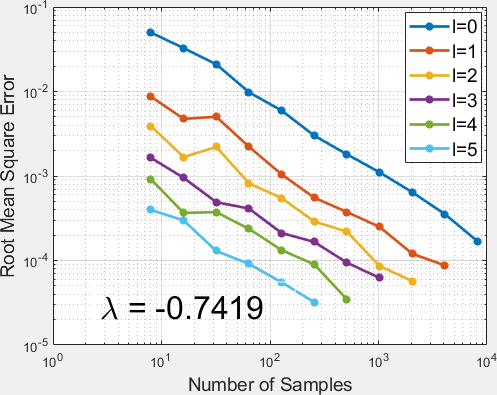} 
    \end{minipage}
    \begin{minipage}{0.33\textwidth}
        \centering
        \includegraphics[width=\textwidth]{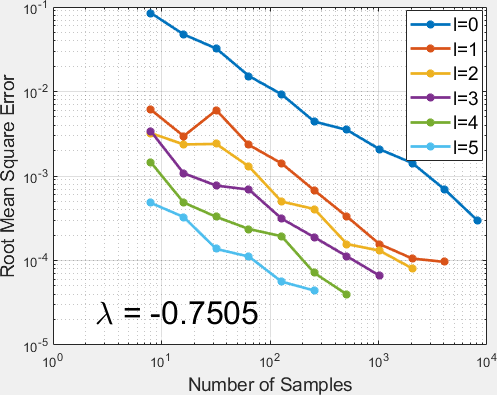}
    \end{minipage}\hfill
    \end{center}
    \caption{Variance test of MLQMC-Sobol for all levels. The variance of the estimator is plotted as a function of the number of samples $N_l$. The convergence rates $\lambda$ are marked on the plot}
    \label{fig:mlqmc_sobol_2_test}
\end{figure}

Next, the results of MLQMC-Sobol', the convergence test (Fig. \ref{fig:mlqmc_sobol_2_test})

\begin{figure}[H]
    \begin{center}
    \hfill
    \begin{minipage}{0.33\textwidth}
        \centering
        \includegraphics[width=\textwidth]{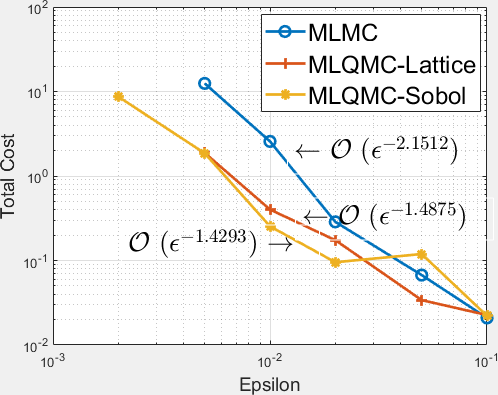} 
    \end{minipage}
    \begin{minipage}{0.33\textwidth}
        \centering
        \includegraphics[width=\textwidth]{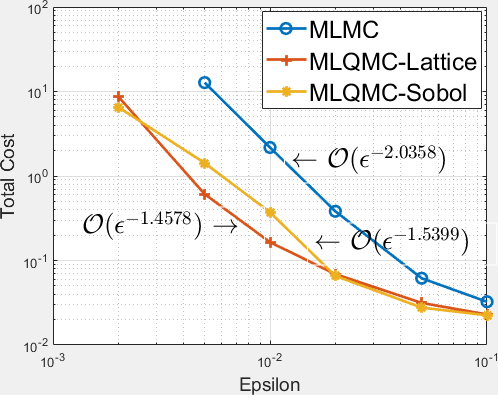} 
    \end{minipage}\hfill
    \null \\
    \begin{minipage}{0.33\textwidth}
        \centering
        \includegraphics[width=\textwidth]{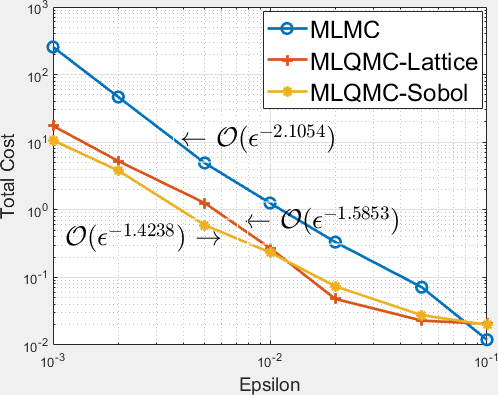} 
    \end{minipage}
    \begin{minipage}{0.33\textwidth}
        \centering
        \includegraphics[width=\textwidth]{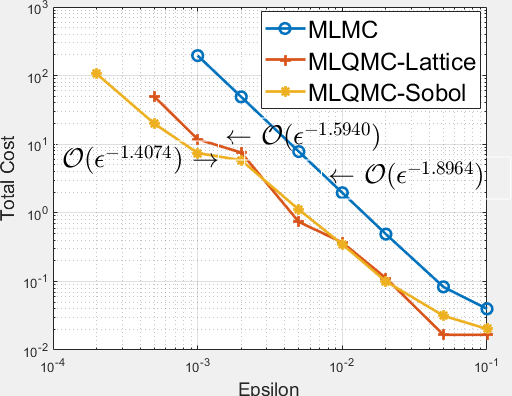}
    \end{minipage}\hfill
    \end{center}
    \caption{Computational complexity of the three methods}
    \label{fig:case2_comparison}
\end{figure}

Finally, we compare the computational cost of the three methods, given $\epsilon$. The results in 4 random fields are presented in Fig. \ref{fig:case2_comparison}. 

\section{Conclusions}\label{Sec4}

In this work, we combined the MLMC with a full multigrid method. We saved the computation for the coarse grid solution on each level without modifying MLMC hierarchy and introducing correlation between different levels. We applied the consistent coarsening approach such that no bias error was introduced in the telescoping formula \eqref{eq:Chap2_Telescoping_Formula}. 

We tested our FMG-MLQMC algorithm on 2-D elliptic PDE with random coefficients for two different types of boundary condition settings and QoIs. The random coefficients were modelled as lognormal random fields with the Mat\'ern covariance function with various parameter settings. We observed that quasi-Monte Carlo approaches have better performance on smoother random fields and problems with more regularity. Also, the comparison of Monte Carlo and quasi-Monte Carlo methods (including Lattice rule and Sobol' sequence) revealed that QMC outperforms MC due to a smaller estimator variance, and Sobol' sequence performs slightly better than Lattice rule.



One future work could be the substitution of the geometric multigrid solver with the algebraic multigrid (AMG) solver. In the AMG scheme, the grids are not associated with physical meshes, rather, the grids are fully determined by the matrix entries algebraically. 

Another work could be to extend the elliptic model to more sophisticated models, such as two-phase porous flow. The efficient sampling and fast simulation of multi-phase subsurface flow under heterogeneous media could produce practical values.

Further work could extend the multilevel model to a multiscale model, and multiscale meshes rather than geometric meshes would be used. In this case, one level could correspond to one scale, and sampling would be performed on each scale. The literature on multiscale modeling can be found in \cite{jenny2003multi,efendiev2009multiscale}. 

\section*{Acknowledgements}
The authors gratefully acknowledge the support from the National Natural Science Foundation of China (Nos. 51874262, 51904031) and the Research Funding from King Abdullah University of Science and Technology (KAUST) through the grants BAS/1/1351-01-01.

\bibliography{References}

\end{document}